\documentclass[twoside]{article}
\usepackage[tbtags]{amsmath}
\usepackage{amsfonts}
\usepackage[mathcal]{eucal}
\usepackage{amssymb}
\usepackage[pdfstartview=FitH,pdfnewwindow=true,pdftitle={The next variational prolongation of the Euclidean space},pdfdisplaydoctitle=true,colorlinks,pagebackref,bookmarksopen,breaklinks]{hyperref}
\newtheorem{prop}{Proposition}
\makeatletter
\def\ps@firstpage{\ps@plain
  \def\@oddfoot{\normalfont\scriptsize \hfil\thepage\hfil}
  \def\@oddhead{\article@logo\hss}}
\def\article@logo{%
  \vbox to\headheight{%
    \@parboxrestore \@logofont
    \noindent
Tensor, N.~S.
    \newline Vol.~68 (2007), no.~1, p.~1--9.
    \par\vss
  }}
\def\@logofont{\fontsize{8}{9.6\p@}\selectfont}
\makeatother
\newbox\stacka\setbox\stacka=\hbox{$\stackrel{\scriptscriptstyle\downarrow}{{\boldsymbol{\mathsf A}}}$}
\edef\stack{\vrule width 0pt depth 0pt height \ht\stacka}
\def\sa{\copy\stacka}

\newcommand{\bmat}[1]{{\boldsymbol#1}}
\newcommand{\sss}[1]{{\scriptscriptstyle#1}}
\newcommand\so{{\scriptscriptstyle 0}}
\newcommand\R{{\mathbb R}}
\newcommand\E{{\mathsf E}}
\newcommand\cL{{\mathcal L}}
\newcommand\cE{{\mathcal E}}
\newcommand\sA{\mathsf A}
\newcommand\sB{\mathsf B}
\newcommand\sC{\mathsf c}
\newcommand\sk{\mathsf k}
\newcommand{\A}{\boldsymbol{\mathsf A}}
\newcommand{\B}{\boldsymbol{\mathsf  B}}
\newcommand{\CC}{\boldsymbol{\mathsf  c}}
\newcommand{\V}{\boldsymbol{\mathsf  v}}
\newcommand{\X}{\boldsymbol{\mathsf  x}}
\newcommand{\K}{\boldsymbol{\mathsf  k}}
\newcommand{\W}{\boldsymbol{\Omega}}

\newcommand{\bP}{\boldsymbol{\pi}}
\newcommand{\Zero}{\boldsymbol{\mathsf  0}}
\newcommand{\bcdot}{\boldsymbol{\cdot}}
\newcommand{\bpartial}{\boldsymbol{\partial}}
\newcommand{\bkey}[1]{\boldsymbol{#1}}
\newcommand{\bs}[1]{\boldsymbol{\mathsf{#1}}}

\begin{document}
\title{THE NEXT VARIATIONAL PROLONGATION OF THE EUCLIDEAN SPACE.}
\author{by Roman \sc Matsyuk}
\date{}
\maketitle
\thispagestyle{firstpage}
\begin{abstract}
The unique third-order invariant variational equation in
three-dimensional (pseudo)Euclidean space is derived.
\end{abstract}
\renewcommand{\thefootnote}{}
\footnotetext{Received December 1, 2005 and in the revized form December 9, 2005.}
\footnotetext{The work was presented in the 8\textsuperscript{th} International Conference of Tensor Society held at Varna, Bulgaria, Aug 22--26, 2005.}
{\bf Introduction.\hspace{.4em}}
The puzzle of what construction should be taken as the local model
of a higher-order Kawaguchi space up to now does not have a common
solution. From the point of view of extremal paths approach one
may assert that for the role of the candidate for the very next
order generalization of the celebrated (pseudo)Riemannian geometry
might be taken a space, the exremals in which do satisfy a
third-order differential equation. In plain geometry the first one
such equation which drops in on one's mind is  that of the
(geodesic) circle, $\frac{d k}{d s}=0$. In 1969 Ukrainian
mathematician Skorobohat'ko suggested building up a geometry in
the Euclidean plain, where geodesics should pass through $n$
arbitrarily given points and therefore be solutions of a
higher-order differential equation of the type $\frac{d^r k}{d
s^r}=0$ for some $r$~\cite{matsyuk:Skoro1, matsyuk:Skoro2}.
Inspired by these ideas I tried to solve the inverse variational
problem for a third order differential equation in
(pseudo)Euclidean space. The fact that one starts from the
(pseudo)Euclidean geometry suggests that the higher-order equation
of geodesic paths one looks for should inherit (pseudo)Euclidian
symmetry. It turns out that in case of three-dimensional space
this problem admits definite solution. Due to the very kind
support of the Organizes of the Conference I take this opportunity
to present the corresponding statements at this talk. The
construction to be proposed here deviates from the notion of
Kawaguchi space in that there does not exist an intrinsically
defined integrand for the variational problem, although the
variational Euler-Poisson equation itself is well defined. On
other hand, to produce a strictly third-order equation, the
integrand should be affine in second order derivatives. This
latter feature relates it to the case of special Kawaguchi space
with $p=1$. Also in three-dimensional space the (vector)
variational equation is necessary degenerate, so one may chose to
prefer that of parameter-indifference generacy, again thus meeting
the terms of Kawaguchi space. I dedicate this special case of
third-order variational space to the name of professor Vitaliy
Skorobohat'ko.

{\bf \S\;1.\hspace{.4em} Preliminary agreements.\hspace{.4em}}
The shortest constructive way to treat the inverse problem of
variational calculus is to introduce the operator of Lagrange
differential $\delta$. In calculable form it was done by Tulczyjew
in~\cite{matsyuk:Tulcz} for the autonomous variational problem and
modified by Kol\'a\v r in~\cite{matsyuk:Kolar} for the case of
non-autonomous one. As our substantial considerations here will
concentrate on (pseudo)Euclidean case, we shall not emphasize
general significance of the notions introduced, but the Reader
will easily understand, what constructions work perfectly well on
general differential manifolds. Let, therefore an
$n+1$-dimensional manifold $M$ be parameterized by local
coordinates $t\equiv x^{\so}$, $x^i$, $i\in\overline{1,n}$ and
consider the space of $r$-th order velocities $T^r
M=J^r(\R,M)(0)$, those being $r$-th order jets with  source zero
from $\R$ to $M$. Let $x^\alpha, u^\alpha, \dot u^\alpha, \ldots
u_{r-1}^\alpha$, $\alpha\in\overline{0,n}$, be local coordinates
in $T^r M$. Germs of one-dimensional submanifolds in $M$ give rise
to another space---that of one-dimensional contact elements in
$M$, denoted by $C^r(M,1)$. This latter space locally is
parameterized by coordinates $t, x^i, v^i, v'^i, \ldots
v_{r-1}^i$. From time to time notations $u_{-\sss1}^\alpha$ and
$v_{-\sss1}^i$ will be used instead of $x^\alpha$ and $x^i$
respectively. The projection

\begin{equation}\label{matsyuk:p}
\wp:\tilde T^r M\to C^r(M,1)
\end{equation} from non-zero velocities space $\tilde T^r M$ to the space of
contact elements is that of quotient projection under the right action of
   the reparametrization group $GL^r(\R)$ on the space $T^r M$. In the third
order this projection is given by the expressions

\begin{equation}\label{matsyuk:p3}
\begin{gathered}
 v^i=\dfrac{u^i}{u^{\so}}\\
 v'^i=\dfrac{\dot u^i}{(u^{\so})^2} - \dfrac{\dot u^{\so}}{(u^{\so})^3}u^i
\\
 v''^i=\dfrac{\ddot u^i}{(u^{\so})^3} - 3\,\dfrac{\dot
 u^{\so}}{(u^{\so})^4}\dot u^i + 3\,\dfrac{(\dot u^{\so})^2}{(u^{\so})^5}u^i
 - \dfrac{\ddot u^{\so}}{(u^{\so})^4}u^i\,.
\end{gathered}
\end{equation}
The generalization of these formulae to arbitrary order of jets
may be found in~\cite{matsyuk:DGA8}. It can be deduced from
general transformation rules for higher order derivatives as
presented, for instance, in~\cite{matsyuk:MKaw}. The contact
elements manifold $C^r(M,1)$ locally is built as the  jet bundle
$J^r(\R,\R^n)$. Any local Lagrange density on $C^r(M,1)$ is therefore best
represented by a semi-basic differential one-form

\begin{equation}\label{matsyuk:UJP11}
\Lambda=L(t, x^i, v^i, \ldots v_{k-1}^i)dt\,.
\end{equation}
The corresponding local Euler-Poisson equations,

\begin{equation}\label{matsyuk:UJP13}
\E_i(t, x^i, v^i, \ldots v_{r-1}^i)=0
\end{equation}
naturally fit in with the conception of a vector differential
one-form

\begin{equation}\label{matsyuk:UJP23}
e=\E_i dx^i\otimes dt\,.
\end{equation}

On the space $T^k M$ one may pose an autonomous variational problem by
introducing a Lagrange function $\cL(x^\alpha, u^\alpha, \ldots
u_{k-1}^\alpha)$, the corresponding Euler-Poisson equations of which,

\begin{equation}\label{matsyuk:calE}
\cE_\alpha(x^\alpha, u^\alpha, \ldots u_{r-1}^\alpha)
\end{equation}
fall
into the shape of globally well defined differential form

\begin{equation}\label{matsyuk:UJP27}
\varepsilon=\cE_\alpha dx^\alpha\,.
\end{equation}
The following assertion is true:

\begin{prop}\label{matsyuk:prop1}
The differential forms $e$ from~(\ref{matsyuk:UJP23}) and

\begin{equation}\label{matsyuk:UJP21}
\varepsilon=-u^i\E_i dx^{\so} + u^{\so}\E_i dx^i
\end{equation}
both satisfy variational criterion simultaneously, if either does.
The corresponding local Lagrangians are related by the formula

\begin{equation}\label{matsyuk:UJP19}
\cL=u^{\so}L\,.
\end{equation}
 \end{prop}
In addition, one observes that the function~(\ref{matsyuk:UJP19})
satisfies the Zermelo conditions, and each such $\cL$ passes to
quotient along the projection~(\ref{matsyuk:p}). A few words on
variational criteria deserve saying then.

{\bf \S\;2.\hspace{.4em} Variational criterion.\hspace{.4em}}
One reason for casting the system of Euler-Poisson equations in the shape of
exterior differential forms is that in the algebra of differential
forms the operator $\delta$ called Lagrange differential may be
introduced. It satisfies $\delta^2=0$, due to what the criterion
of the existence of a local Lagrange function for, say, the system
of equations~(\ref{matsyuk:UJP13}) is expressed as $\delta e=0$
for $e$ in~(\ref{matsyuk:UJP23}).

Consider the graded algebra of differential forms on the space
$J^r(\R,Q)\approx\R\times T^r Q$ of jets from $\R$ to arbitrary
manifold $Q$. Let us recall that an operator $D$ is called a
derivation of degree $q$ if for any differential form $\varpi$ of
degree $p$ and any other differential form $\omega$ the
differential form $D\varpi$ is of degree $p+q$ and the Leibniz rule
$D(\varpi\wedge\omega)=D\varpi\wedge\omega + (-1)^{pq}\varpi\wedge
D\omega$ holds. Let us recall some familiar operators acting on
forms. The operator of vertical differential $d_v$ is first
defined on the ring of functions as $d_v f=\sum_i\frac{\partial
f}{\partial x^i}dx^i + \sum_{s=0}^{r-1}\sum_i\frac{\partial
f}{\partial v_s^i}dv_s^i$, $\{x^i\}\in Q$, and then extended as a derivation of
degree 1 by means of the coboundary property $d_v^2=0$. The total
derivative $D_t$ is also first defined on the ring of functions as
$D_t f=\frac{\partial f}{\partial t} + \sum_i v^i\frac{\partial
f}{\partial x^i} + \sum_{s=0}^{r-1}\sum_i v_{s+1}^i\frac{\partial
f}{\partial v_s^i}$, and then extended as a derivative of degree
zero by means of the commutation relation $D_td_v=d_vD_t$.
Following Tulczyjew, we need one more derivation of degree zero,
denoted here by $\iota$, and defined by its action on functions
and forms as $\iota f=0$, $\iota dx^i=0$, $\iota
dv_s^i=(s+1)dv_{s-1}^i$, $s\in\overline{0,r-1}$.  Let us denote by
$\iota^{\so}$ the operator of evaluating the degree of a
differential form and by $D^s$ the iterated $D$.  The Lagrange
$\delta$ is first introduced by its action in the algebra of differential forms
on $T^r Q$, eventually with coefficients depending on the time $t\in\R$,

\[
\delta=\sum_{s=0}^r\frac{(-1)^s}{s!}D_t^s\iota^s d_v\,,
\]
and afterwards trivially extended to the graded module of semi-basic with
respect to $\R$ differential forms on $J^r(\R,Q)$ (actually one-forms) with
coefficients in the bundle of graded algebras $\wedge T^*(T^r Q)\to T^r Q$
by means of the prescriptions:

\[
\delta(\omega\otimes dt)= \delta(\omega)\otimes dt \,.
\]
The property $\delta^2=0$ holds. One may apply either the notion
of the (above defined time-extended) Lagrange differential to forms on the jet space
$J^r(\R,\R^n)$, setting $Q=\R^n$, or the notion of the ``truncated'' time-independent Lagrange
differential to the forms both on the manifold $T^r M$ as well as on the
manifold $T^r(\R\times\R^n)$ setting $Q=M$ and $Q=\R\times\R^n$
respectively. Thus locally the notion of the Lagrange differential
is applicable to both sides of the projection~(\ref{matsyuk:p}),
whereas globally it is well defined on the left hand side solely.
In each case the differential  forms (\ref{matsyuk:UJP23})
and~(\ref{matsyuk:UJP27}) that represent the Euler-Poisson
equations are in fact semi-basic also with respect to $Q$. In terms
of the operators introduced above this means that $\iota e$ and
$\iota \varepsilon$ both are zero.

\begin{prop}
Considering formulae~(\ref{matsyuk:UJP11}), (\ref{matsyuk:UJP23}),
(\ref{matsyuk:UJP21}) and (\ref{matsyuk:UJP19}), if $e=\delta\Lambda$, then
$\varepsilon=\delta\cL$. The variational criterion for~(\ref{matsyuk:UJP23})
consists in $\delta e=0$ and is equal to $\delta\varepsilon=0$.
 \end{prop}

The criterion $\delta e=0$ now can be expressed in coordinates.
After some permutations of indices and some interchanges in the
order of sequential sums one gets in a way similar to that
of~\cite{matsyuk:Lawruk}

\[
\delta e=\sum_{s=0}^r
\left(\dfrac{\partial\E_i}{\partial
v_{s-1}^j}-\sum_{k=s}^r(-1)^k\dfrac{k!}{(k-s)!s!}D_t^{k-s}\dfrac{\partial\E_j}{\partial
v_{k-1}^i}\right)dv_{s-1}^j\wedge dx^i\,,
\] from where the following system of partial differential equations
follows:  \begin{subequations}\label{matsyuk:2}
\begin{align}
\label{matsyuk:2.1}
 \dfrac{\partial\E_i}{\partial x^j}
 - \dfrac{\partial\E_j}{\partial x^i}
+\sum_{k=0}^r(-1)^k D_t^k
\left(
\dfrac{\partial\E_i}{\partial v_{k-1}^j}
-\dfrac{\partial\E_j}{\partial v_{k-1}^i}
\right) &=0\,;&
 \\
\label{matsyuk:2.2}
\dfrac{\partial\E_i}{\partial
v_{s-1}^j}-\sum_{k=s}^r(-1)^k\dfrac{k!}{(k-s)!s!}D_t^{k-s}\dfrac{\partial\E_j}{\partial
v_{k-1}^i}&=0& 1\leqslant s\leqslant r\,.
\end{align}
\end{subequations}

The above system of equation is equivalent to the following one (obtained from~(\ref{matsyuk:2.2})
by extending the range of $s$ to
include $s=0$):

\begin{align}\label{matsyuk:criterion}
\dfrac{\partial\E_i}{\partial
v_{s-1}^j}-\sum_{k=s}^r(-1)^k\dfrac{k!}{(k-s)!s!}D_t^{k-s}\dfrac{\partial\E_j}{\partial
v_{k-1}^i}&=0& 0\leqslant s\leqslant r\,.
\end{align}

{\it Proof.\hspace{.4em}}
The antisymmetrization of~(\ref{matsyuk:criterion}) at $s=0$ produces the equation~(\ref{matsyuk:2.1}).
On the contrary, in equation~(\ref{matsyuk:2.1}) separate the summand with
$k=0$:

 \begin{equation}\label{matsyuk:2.1.3}
 2\,\dfrac{\partial\E_i}{\partial x^j}
 - 2\,\dfrac{\partial\E_j}{\partial x^i}
+\sum_{k=1}^r(-1)^k D_t^k
\dfrac{\partial\E_i}{\partial v_{k-1}^j}
-\sum_{k=1}^r(-1)^k D_t^k\dfrac{\partial\E_j}{\partial v_{k-1}^i}
 =0\,.
\end{equation}
Under the first sum sign substitute
$\dfrac{\partial\E_i}{\partial v_{k-1}^j}$ from
equation~(\ref{matsyuk:2.2}):

\[
\sum_{k=1}^r(-1)^k D_t^k \dfrac{\partial\E_i}{\partial v_{k-1}^j}
=\sum_{k=1}^r(-1)^k D_t^k
\sum_{s=k}^r(-1)^s\dfrac{s!}{(s-k)!k!}D_t^{s-k}\dfrac{\partial\E_j}{\partial
v_{s-1}^i}\,.
\]
Interchange the summation order:
$\sum_{k=1}^r\sum_{s=k}^r=\sum_{\substack{ s,k=1
\\ s\geqslant k}}^r=\sum_{s=1}^r\sum_{k=1}^s$. Calculate the sum
over $k$:

\[
\sum_{k=1}^s(-1)^k\dfrac{s!}{(s-k)!k!}=\sum_{k=0}^s(-1)^k \binom s k -
\binom s 0 =0-1=-1\,.
\]
Ultimately equation~(\ref{matsyuk:2.1})
becomes \[ 2\,\dfrac{\partial\E_i}{\partial x^j} -
2\,\dfrac{\partial\E_j}{\partial x^i} -\sum_{k=1}^r(-1)^k D_t^k
 \dfrac{\partial\E_j}{\partial v_{k-1}^i} -\sum_{k=1}^r(-1)^k
D_t^k\dfrac{\partial\E_j}{\partial v_{k-1}^i} =0\,, \] which
coincides with doubled equation~(\ref{matsyuk:criterion}) at $s=0$.  The
criterion~(\ref{matsyuk:criterion}) has been obtained by different authors. The Reader may
consult the book~\cite{matsyuk:OKrup} by Olga Krupkov\'a for a
recent review.

Let us focus on third order variational equations.
It is obvious that the Euler-Poisson expressions
are of affine type in the highest derivatives.
We utilize some familiar vector notations: the lower dot symbol
will denote the contraction between a row-array and the subsequent
column-array and sometimes also will stand for the matrix multiplication
between a matrix and the subsequent column-array.  From the system of partial differential equations~(\ref{matsyuk:criterion})
it is possible to deduce that
the most general form of the Euler-Poisson equation of the third order reads:

\begin{equation}\label{matsyuk:hamspin5}
\A{\,\bkey.\,}\V^{{{\prime}}{{\prime}}}{\, +
\,}(\V^{{\prime}}{\!\bkey.\,}{\bpartial}_{\V})\,
\A{\,\bkey.\,}\V^{{\prime}}{\, + \,}\B{\,\bkey.\,}\V^{{\prime}}{\, +
\,}\CC\, = \,\Zero\,,
\end{equation}
where the skew-symmetric matrix
$\A$, the symmetric matrix $\B$, and a column $\CC$ all depend on $t$, $x^i$,
and $v^i$ and satisfy the following system of partial differential
equations:

\begin{subequations}\label{matsyuk:hamspin6}
{\renewcommand{\arraystretch}{1.7}
\begin{gather}
\label{matsyuk:hamspin6.1}        \partial_{_{_{_{{{v}}}}}}{\!}_{[i}{}{{\sA}}_{jl]}=0
\\ \label{matsyuk:hamspin6.2}        2\,{{\sB}}_{[ij]}-3\,{\bf D_{_{\bmat 1}}}{\kern.01667em}{{\sA}}_{ij}=0
\\ \label{matsyuk:hamspin6.3}     2\,\partial_{_{_{_{{{v}}}}}}{\!}_{[i}{}{{\sB}}_{j]\,l}
               -4\,\partial_{_{_{_{{{x}}}}}}{\!}_{[i}{}{{\sA}}_{j]\,l}
               +{\partial_{_{_{_{{{x}}}}}}{\!}_{l}}{\,}{{\sA}}_{ij}
   +2\,{\bf D_{_{\bmat 1}}}{\kern.01667em}{\partial_{_{_{_{{{v}}}}}}{\!}_{l}}{\,}{{\sA}}_{ij}=0
\\ \label{matsyuk:hamspin6.4}       {\partial_{_{_{_{{{v}}}}}}{\!}_{(i}}{}{{\sC}}_{j)}
               -{\bf D_{_{\bmat 1}}}{\kern.01667em}{{\sB}}_{(ij)}=0
\\ \label{matsyuk:hamspin6.5}       2\,{\partial_{_{_{_{{{v}}}}}}{\!}_{l}}{\,}\partial_{_{_{_{{{v}}}}}}{\!}_{[i}{}{{\sC}}_{j]}
           -4\,\partial_{_{_{_{{{x}}}}}}{\!}_{[i}{}{{\sB}}_{j]\,l}
           +{{\bf D_{_{\bmat 1}}}}^{2}{\,}{\partial_{_{_{_{{{v}}}}}}{\!}_{l}}{\,}{{\sA}}_{ij}
 +6\,{\bf D_{_{\bmat 1}}}{\kern.0334em}\partial_{_{_{_{{{x}}}}}}{\!}_{[i}{}{{\sA}}_{jl]}=0
\\ \label{matsyuk:hamspin6.6}
4\,\partial_{_{_{_{{{x}}}}}}{\!}_{[i}{}{{\sC}}_{j]} -2\,{\bf D_{_{\bmat
           1}}}{\kern.0334em}\partial_{_{_{_{{{v}}}}}}{\!}_{[i}{}{{\sC}}_{j]} -{{\bf D_{_{\bmat 1}}}}^{3}{\,}{{\sA}}_{ij}=0\,.
\end{gather}}
\end{subequations}
Here the differential operator $\bf D_{_{\bmat 1}}$ is the lowest
order truncated operator of total derivative $D_t$,
\[
{\bf D_{_{\bmat 1}}}=\partial_{t}{\,+\,}\V{\,\bkey.\,}{\bpartial}_{\X}\,.
\]

Alongside with the differential form~(\ref{matsyuk:UJP23})  it is
convenient to introduce the so-called Lepagian equivalent to it,
whose coefficients do not depend on third-order derivatives:

\begin{equation}\label{matsyuk:UJP31}
\begin{split}
\epsilon=\sA_{ij} dx^i\otimes dv^{\prime}{}^j + &\, \sk_i dx^i \otimes dt, \quad\mathrm{where}
\\ &\K = (\V^{{\prime}}{\!\bkey.\,}{\bpartial}_{\V})\,
\A\,\bkey.\,\V^{{\prime}}{\, + \,}\B\,\bkey.\,\V^{{\prime}}{\, +
\,}\CC\,.
\end{split}
\end{equation}
This vector-valued differential one-form (taking values in
$T^*\R^n$) may be thought of as an interpretation of the Lepagian
form, alternative to that considered in~\cite{matsyuk:OKrup}.

Since we are interested in holonomic local curves in~$C^3(M,1)$, it is a
common point that the vector-valued differential
one-forms~(\ref{matsyuk:UJP31}) and (\ref{matsyuk:UJP23}) are treated as equal
with respect to the contact module on $J^3(\R,\R^n)$:

\[
\epsilon-e =\sA_{ij}dx^i \otimes   \theta_{\bmat3}^j\,,
\]
where the vector-valued contact one-forms

\begin{equation}
\label{matsyuk:hamspin8}
\boldsymbol{\theta_1}=d\X-\V dt\,,\quad\boldsymbol{\theta_2}= d \V-\V' dt\,,
\quad \boldsymbol{\theta_3}=d\V'-\V'' dt
\end{equation}
generate the contact module on $J^3(\R,\R^n)$.

{\bf \S\;3.\hspace{.4em} Euclidean symmetry.\hspace{.4em}}
Since we simultaneously consider both the true Euclidean and the
pseudo-Euclidean cases, let us fix some notations. By $\eta$ the
sign $+$ or $-$ of the component $g_{\so\so}$ of the canonical
diagonal metric tensor will be denoted. Centered dot will mean
scalar product between matrices which represent tensors or between
arrays which represent vectors---with respect to the
(pseudo)Euclidean canonical metric tensor. Thus the scalar product
is merely the contraction that involves the metric tensor. The
infinitesimal generator $X$ of the (pseudo)Euclidean
transformation in three-dimensional space may be parametrized by
means of a skew-symmetric matrix $\W$ and some vector $\bP$:

\eqnarray
 X&=&
-(\bP\bcdot{\X})\,{\partial_t}+\eta\,t\,\bP\,\bkey.\,\bpartial_{\X}
+\W\bcdot({\X}\wedge\bpartial_{\X})
\nonumber\\ &&{}+\eta\,\bP\,\bkey.\,\bpartial_{\V}
+(\bP\bcdot \V)\,\V\,\bkey.\,\bpartial_{\V}+\W\bcdot(\V\wedge\bpartial_{\V})\nonumber
\\ \nonumber &&{}+2\,(\bP\bcdot \V)\,\V'\bkey.\,\bpartial{_{\V'}}+(\bP\bcdot \V')\,\V\,\bmat.\,\bpartial{_{\V'}}
+\W\bcdot(\V'\wedge\bpartial{_{\V'}})\,.
\endeqnarray

It is possible to cast the idea of the symmetry of the equation~(\ref{matsyuk:hamspin5})
into the framework of exterior differential system invariance concept.
The system to handle is generated by the vector-valued Phaff  form $\bmat\epsilon$ from~(\ref{matsyuk:UJP31}) and
the contact vector-valued differential forms $\boldsymbol{\theta_1}$ and
$\boldsymbol{\theta_2}$ from~(\ref{matsyuk:hamspin8}).
Let $X(\bmat\epsilon)$ denote the Lie derivative of the vector-valued differential form
$\bmat\epsilon$ along the vector field $X$. The invariance condition consists in that
there may be found some matrices ${\bf\Phi}$, ${\bf\Xi}$, and ${\bf\Pi}$
depending on $\V$ and $\V'$ such that
\begin{equation}\label{matsyuk:hamspin9}
X(\bmat\epsilon)={\bf\Phi}\,\bmat.\,\bmat\epsilon+{\bf\Xi}\,.\,(d\X-\V d t)
+{\bf\Pi}\,\bmat.\,(d\V-\V' d t).
\end{equation}
We also assert that $\A$ and $\K$ in~(\ref{matsyuk:UJP31}) do not depend neither
on $t$ nor on $\X$.

The identity~(\ref{matsyuk:hamspin9}) splits into more identities, obtained by evaluating the
coefficients of the differentials $d t$, $ d\X$,
$d\V$, and $d\V'$ independently:
\begin{eqnarray}
\nonumber
&{\big(\bP\,\bmat.\,\bpartial_{\V}+(\bP\bcdot \V)\,\V\,\bmat.\,\bpartial_{\V}+\W\bcdot(\V\wedge \bpartial_{\V})\big)\,\A}
\hspace{2.9cm}
&\\ \label{matsyuk:81}&{{}+ 2\,(\bP\bcdot \V)\A+(\A\V)\otimes\bP-\A\W={\bf\Phi}\A}\,;&
\\ \label{matsyuk:82}& 2\,(\A\V')\otimes\bP+(\bP\bcdot \V')\,\A\,={\bf\Pi}\,;&
\\ \label{matsyuk:83}&-{\K}\otimes\bP={\bf\Xi}\,;&
\\ \label{matsyuk:84}&  X({\K})={\bf\Phi}{\K}-{\bf\Xi}\V-{\bf\Pi}\V'\, .&
\end{eqnarray}
In the above the `$\otimes$' symbol means the tensor (sometimes named as
`direct') product of matrices; the associative matrix multiplication
is represented by joint writing.

A skew-symmetric two-by-two matrix always has the inverse, so the
`Lagrange multipliers' ${\bf\Phi}$, ${\bf\Xi}$, and ${\bf\Pi}$ may explicitly be
defined from the equations (\ref{matsyuk:81}--\ref{matsyuk:83}) and then
substituted into~(\ref{matsyuk:84}).  Subsequently, the equation
(\ref{matsyuk:84}) splits into the following identities by the powers of the
variable $\V'$ and by the parameters $\W$ and $\bP$ (take notice of the
derivative matrix $\A^{\prime}=(\V'\bmat.\,\bpartial_{\V})\,\A$; also the
vertical arrow sign points to the very last factor to which the aforegoing
differential operator still applies):

\begin{eqnarray}
\nonumber
     \lefteqn{\bigl(\W\bcdot(\V\wedge \bpartial_{\V})\bigr)\,\A^{\prime}\V'+\bigl(\W\bcdot(\V'\wedge \bpartial_{\V})\bigr)\,\A\V'
     -(\V'\bmat.\,\bpartial_{\V})\,\A\W\V'}
\\   \label{matsyuk:95}
         &{{}=\bigl(\W\bcdot(\V\wedge \bpartial_{\V})\bigr)\,\sa \A^{-1}\A^{\prime}\V'-\A\W\A^{-1}\A^{\prime}\V'}\,;
&\\ &\label{matsyuk:96}\stack\bigl(\W\bcdot(\V\wedge \bpartial_{\V})\bigr)\,\B-\B\W
     =\bigl(\W\bcdot(\V\wedge \bpartial_{\V})\bigr)\,\sa \A^{-1}\B-\A\W\A^{-1}\B\,;
&\\ &\label{matsyuk:97}\stack\bigl(\W\bcdot(\V\wedge \bpartial_{\V})\bigr)\,\CC
     =\bigl(\W\bcdot(\V\wedge \bpartial_{\V})\bigr)\,\sa \A^{-1}\CC-\A\W\A^{-1}\CC\,;
&\\ \nonumber
     \lefteqn{\stack\bigl(\bP\,\bmat.\,\bpartial_{\V}+(\bP\bcdot \V)\,\V\,\bmat.\,\bpartial_{\V}\bigr)\,\A^{\prime}\V'
     +(\bP\bcdot \V)\,\A^{\prime}\V'+(\bP\bcdot \V')\,(\V\,\bmat.\,\bpartial_{\V})\,\A\V'+(\bP\bcdot \V')\,\A^{\prime}\V}
\\ &   \label{matsyuk:98}{}=\bigl(\bP\,\bmat.\,\bpartial_{\V}+(\bP\bcdot \V)\,\V\,\bmat.\,\bpartial_{\V}\bigr)\,\sa \A^{-1}\A^{\prime}\V'
          +(\bP \A^{-1}\A^{\prime}\V')\,\A\V- 3 \,(\bP\bcdot \V')\,\A\V'\,;
&\\ \nonumber
     \lefteqn{\stack\bigl(\bP\,\bmat.\,\bpartial_{\V}
     +(\bP\bcdot \V)\,\V\,\bmat.\,\bpartial_{\V}\bigr)\,\B+(\B\V)\otimes\bP}
\\ &\label{matsyuk:99}
        {}=\bigl(\bP\bcdot\bpartial_{\V}+(\bP\bcdot \V)\,\V\,\bmat.\,\bpartial_{\V}\bigr)\,\sa \A^{-1}\B
        +(\A\V)\otimes\bP \A^{-1}\B+(\bP\bcdot \V)\B\,;
&\\ \nonumber
     \lefteqn{\stack\bigl(\bP\,\bmat.\,\bpartial_{\V}+(\bP\bcdot \V)\,\V\,\bmat.\,\bpartial_{\V}\bigr)\,\CC}
\\ &\label{matsyuk:90}{}=\bigl(\bP\,\bmat.\,\bpartial_{\V}+(\bP\bcdot \V)\,\V\,\bmat.\,\bpartial_{\V}\bigr)\,\sa \A^{-1}\CC
     + 3 \,(\bP\bcdot \V)\,\CC+(\bP \A^{-1}\CC)\,\A\V\,.
\hspace{-1cm}&\hspace{1cm}\end{eqnarray}

Straightforward but cumbersome routine calculations accompanying
the simultaneous solving of the partial differential equations
(\ref{matsyuk:95}) and (\ref{matsyuk:98}) with respect to the
unknown function ${\sA}_{12}$ produce the unique output of
\[
{\sA}_{12}=\frac{\rm const}{(1+{v}_{\sss1}{v}^{\sss1}+{v}_{\sss2}{v}^{\sss2})^{3/2}\strut}\,.
\]

We remind that the system of the equations \{(\ref{matsyuk:95})--(\ref{matsyuk:90})\} and the
system~(\ref{matsyuk:hamspin6}) must be solved simultaneously. Thus, the
equation (\ref{matsyuk:hamspin6.1}) becomes trivial now.

Under the assumption of $\B$ being a symmetric matrix (see
(\ref{matsyuk:hamspin6.2})),
the solution of the equations \{(\ref{matsyuk:96}), (\ref{matsyuk:99})\} is:
\[
{\sB}_{ij}={\rm const}\cdot(1+\V\bcdot \V))^{-3/2}\bigl({v}_{i}{v}_{j}
-(1+\V\bcdot \V)\,{g}_{ij}\bigr)\,.
\]

This automatically satisfies  the equation (\ref{matsyuk:hamspin6.3}) too.
In what concerns the subsystem \{(\ref{matsyuk:97}), (\ref{matsyuk:90})\},
only the trivial solution
$\CC={\Zero}$ exists.

We are ready now to formulate the summary of the above  development in terms of a proposition:
\begin{prop}
The invariant parameter-indifferent Euler-Poisson equation in
three-dimensional (pseudo)Euclidean space is:
\begin{equation}
\label{matsyuk:hamspin20}
-\frac{\ast\V''\strut}{(1+\V\bcdot\V)^{3/2}\strut}+3\,\frac{\ast\V'\strut}{(1+\V\bcdot\V)^{5/2}}\,(\V\bcdot\V')
+\frac{\mu\strut}{(1+\V\bcdot\V)^{3/2}}\,\bigl((1+\V\bcdot\V)\,\V'-(\V'\bcdot\V)\,\V\bigr)={\Zero}\,.
\end{equation}
\end{prop}

The arbitrary constant $\mu$ serves to parameterize the set of all the
variational equations (\ref{matsyuk:hamspin20}). The definition of
the `star operator' is common. Thus, $\,\ast1=\bs
e_{(\sss1)}\wedge\,\bs e_{(\sss2)}$, whereas $\,\ast\,(\bs
e_{(\sss1)}\wedge\,\bs e_{(\sss2)})=1\;$ if the
(pseudo)\-\kern.5pt orthonormal frame $\,\{\bs e_{(\sss1)}\,,\ \bs
e_{(\sss2)}\}\,$ carries the positive orientation; also
$\,(\ast{\bs w})_{i}=\varepsilon_{ji}{w}^{j}\,$ for a
two-dimensional vector $\,\bs w$.

I know two different
($j=1,2$) Lagrange functions which produce the
equation~(\ref{matsyuk:hamspin20}),
\begin{equation}
\label{matsyuk:hamspin31} L_{(j)}=\frac{\ast(\V^{\prime}\wedge\bs
e_{(j)})}{(1+\V\bcdot\V)^{1/2}(1+{g}_{jj}\|\V\wedge\bs e_{(j)}\|^{2})}
     {v}^{j}-\mu\, (1+\V\bcdot\V)^{1/2}\,.
\end{equation}
These differ by the total time derivative:
\[
L_{(2)}-L_{(1)}=\frac{d}{dt}\arctan\frac{{v}^{\sss1}{v}^{\sss2}}{\sqrt{1+{v}_{j}{v}^{j}}}\,.
\]

{\it Remark~1.\hspace{.4em}} Equation~(\ref{matsyuk:hamspin20}) describes
helices with second curvature equal to $\|\mu\|$.

{\it Remark~2.\hspace{.4em}} The point
symmetries of the equation~(\ref{matsyuk:hamspin20}) are exhausted
by (pseudo)Euclidean transformations if $\mu\ne 0$. Otherwise they
precisely consist of conformal ones~\cite{matsyuk:theses}.

{\it Remark~3.\hspace{.4em}} There does not exist an invariant
affine second-order Lagrange function in (pseudo)Euclidean space of dimension greater
than $2$ (strictly speaking, this was proved for the signature not equal $2$)~\cite{matsyuk:DAN}.

With the Proposition(\ref{matsyuk:prop1}) in hand and applying formula~(\ref{matsyuk:p3}) it is not
difficult to put down the ``homogeneous''
counterpart~(\ref{matsyuk:calE}) of
equation~(\ref{matsyuk:hamspin20}). It reads:

\begin{equation}\label{matsyuk:UJP36}
-\dfrac{{\bf\ddot{\bmat u}}\times \bmat u\strut}{\|\bmat u\|^3\strut}
+3\,\dfrac{{\bf\dot{\bmat u}}
\times \bmat u\strut}{\|\bmat u\|^5}\,({\bf\dot{\bmat u}}\cdot \bmat u)
-\dfrac{\mu\strut}{\|\bmat u\|^3}\,\left[ (\bmat u\bcdot \bmat u)\,{\bf\dot{\bmat u}} -
({\bf\dot{\bmat u}}\cdot \bmat u)\,\bmat u  \right] = \Zero\,.
\end{equation}

Furthermore, by same means of~(\ref{matsyuk:UJP19}) one may deduce
 a general formula for the family $(\beta=0,1,2)$ of the Lagrange functions
which produce the right hand side of~(\ref{matsyuk:UJP36}):

\begin{equation}\label{matsyuk:DGA95_29}
\cL_{(\beta)}
\;=\;\dfrac{u^{\beta}[{\bf\dot{\bmat u}},\bmat u,\bmat e_{(\beta)}]}
{\|\bmat u\|\,\|\bmat u\times\bmat e_{(\beta)}\|^{2}}
\;-\;m\,\|\bmat u\|
\;+\;{\bf\dot{\bmat u}}\,\bmat.\,\boldsymbol\partial_{\bmat u}\;\phi
\;+\;\bmat a\,\bmat.\,\bmat u\,,
\end{equation}
where an arbitrary row vector $\bkey a$ is constant and a function $\phi$ depending on the variable
$\bkey u$ is subject to the constraint $\bkey u\,\bkey.\,\boldsymbol\partial_{\bkey u}\,\phi\;=\;0$.
Recall also the notation $[\;,\;,\;]$ for the parallelepipedal product of three vectors.
The vector $\bkey e_{(\beta)}$ denotes the $\beta$-th component of the
(pseudo)Euclidean frame.
{\it Each $\cL_{(\beta)}$ fits in.}

\medskip

The problem of finding invariant
variational equations in some special cases,
discussed in this talk, might have been formulated in still more
recent framework of invariant variational bicomplexes (cf. for
example~\cite{matsyuk:Kohan}). Unfortunately, the threshold of the
non-existence of invariant Lagrangian functions diminishes the
effectiveness of the corresponding machinery, which from the very
beginning suggests the invariance of the full bicomplex.
Similar difficulties arise when one starts to apply notions developed for
Kawaguchi spaces. For example, the metric `tensor' calculated from the
Lagrange function~(\ref{matsyuk:DGA95_29}) does not designate any geometric
object. Only quantities, built of the invariant momentum

\[
\mathcal P\stackrel{\rm def}=\dfrac{\partial \cL}{\partial \bmat u}
-\left(\dfrac{\partial \cL}{\partial \bf\dot{\bmat
u}}\right)^{\displaystyle\bcdot} = \dfrac{{\bf\dot{\bmat u}}\times \bmat
u\strut}{\|\bmat u\|^3\strut} +\mu\,\dfrac{\bmat u\strut}{\strut\|\bmat u\|}
\] would play any significant role in a generally covariant theory.
Several such quantities were introduced in chapter~2 of
paper~\cite{matsyuk:AKaw}.

\hfill
\begin{minipage}{15em}\begin{center}
Institute for Applied Problems\\in Mechanics and
Mathematics\\15~Dudayev~St.\\290005 L\kern-1pt'viv, Ukraine
\\E-mail: matsyuk@lms.lviv.ua, romko.b.m@gmail.com
\end{center}\end{minipage}
\rightline{\url{http://www.iapmm.lviv.ua/12/eng/files/st_files/matsyuk.htm}}
\end{document}